\documentclass{article}
\usepackage{graphicx}

\usepackage{amsfonts,amsmath,amsthm}
\newtheorem{definition}{Definition}[section]
\newtheorem{theorem}{Theorem}[section]

\newcommand{\A}{\mathcal F}
\newcommand{\IR}{\mathbb R}

\newcommand{\Cov}{{\rm Cov\,}}

\newcounter{romva}

\newcounter{alphva}

\usepackage{amsmath,amsthm,amsfonts}

\title{IFSM representation of Brownian motion with applications to simulation}

\author{Stefano Maria Iacus\footnote{\texttt{stefano.iacus@unimi.it}, {\it corresponding author.}}, Davide La
Torre\footnote{\texttt{davide.latorre@unimi.it}}\\ Department of Economics, Business and Statistics\\ University of
Milan, Via Conservatorio, 7, I-20122 Milan - Italy\\}

\begin{document}
\maketitle

\begin{abstract}
Several methods  are currently available to simulate paths of the Brownian motion.
In particular, paths of the BM can be simulated using the properties of the increments of the process like in the Euler scheme, or as the limit of a random walk or via $L^2$ decomposition like the Kac-Siegert/Karnounen-Loeve series. In this paper we first propose a IFSM (Iterated Function Systems with Maps) operator whose fixed point is the trajectory of the BM. We then use this representation of the process to simulate its trajectories. The resulting simulated trajectories are self-affine, continuous and fractal by construction. This fact produces more realistic trajectories than other schemes in the sense that their geometry is closer to the one of the true BM's trajectories. The IFSM trajectory of the BM can then be used to generate more realistic solutions of stochastic differential equations.

\end{abstract}

{\bf AMS Subject Classification:}


\section{Introduction}

In this paper we show how to solve the inverse problem
for IFSM in the case of trajectories of stochastic processes in $L^2$.
The method is based on the solution of the inverse problem for IFSM
due to Forte and Vrscay \cite{Forte}.
This is an extension of classical IFS methods which
can be used for approximating a given
element of $L^2(H)$ thus in particular trajectories of stochastic processes on this space.
The final goal of this approach is simulation.
Indeed, several methods  are currently available to simulate paths of stochastic processes and in particular of the Brownian motion.
Paths of the BM can be simulated using the properties of the increments of the process like in the Euler scheme \cite{KPS00}, or as the limit of a random walk or via $L^2$ decomposition like the Kac-Siegert/Karnounen-Loeve series \cite{kac}. In this paper we first propose a IFSM (Iterated Function Systems with Maps) operator whose fixed point is the trajectory of the BM. We then use this representation of the process to simulate its trajectories. The resulting simulated trajectories are self-affine, continuous and fractal by construction. This fact produces more realistic trajectories than other schemes in the sense that their geometry is closer to the one of the true BM's trajectories. The IFSM trajectory of the BM can then be used to generate more realistic solutions of stochastic differential equations.
The paper is organized as follows: Section \ref{sec:ifsm} recalls the theory of IFSM on $L^2$, Section \ref{sec:l2} recalls some details for stochastic processes with trajectories in $L^2$ and the link with the IFSM theory. Section \ref{sec:bm} presents the application of the IFSM theory to the problem of simulation with particular attention to the case of the Brownian motion.

\section{IFS with Maps (IFSM) on $L^2(H)$}\label{sec:ifsm}

The basic idea of Iterated Function Systems (IFS) can be
traced back to some historical papers but the use of such systems
to construct fractals and other similar sets was first described by Hutchinson (1981).
The fundamental result on which the IFS method is based is Banach
theorem. 
The mathematical context is the following: given $y$ in a complete metric space
$(Y,d)$, find a contractive operator $T:Y\to Y$ that admits a
unique fixed point $y^*\in Y$ such that $d_Y(y,y^*)$ is
small enough. In fact if one is able to solve {\it the inverse problem}
with arbitrary precision, it is possible to identify $y$ with the
operator $T$ which has it as fixed point.
The fundamental theorems on which the IFS method is based on
are the following:

\begin{theorem}(Banach Theorem)
Let $(Y,d_Y)$ be a complete metric space; suppose there exists a 
mapping $T:Y\to Y$ such that
$$
d_Y(T(x),T(y)\le c d_Y(x,y)
$$
for all $x,y\in Y$ and some $c\in [0,1)$.
$c$ is said to be the contractivity factor of $T$.
Then there exits
a unique $y^*\in Y$ such that $T(y^*)=y^*$ and for any $y\in Y$
we have $d(T^n(y),y^*)\to 0$ when $n\to +\infty$.
\end{theorem}

\begin{theorem}(Collage Theorem)
Let $(Y,d_Y)$ be a complete metric space. Given $y\in Y$ suppose
that there exists a contractive map $T$ with contractivity factor
$c\in [0,1)$ such that $d_Y(y,T(y))<\epsilon$. If $y^*$ is
the fixed point of $T$ then $d_Y(y,y^*)\le \frac{\epsilon}{1-c}$.
\end{theorem}

\begin{theorem} 
Let $(Y,d_Y)$ be a complete metric space and $T_1,T_2$ be two contractive
mappings with fixed points $y^*_1$ and $y^*_2$. Then
$$
d_Y(y^*_1,y^*_2)\le \frac{1}{1-c_1} d_{Y,\sup}(T_1,T_2)
$$
where
$$
d_{Y,\sup}(T_1,T_2)=\sup_{x\in Y} d(T_1(x),T_2(x))
$$
and $c_1$ is the contractivity factor of $T_1$.
\end{theorem}

We are going to use a particular class of IFS operators, known as IFSM (IFS with Maps), introduced by Forte and Vrscay in 1994. Let $H = [0,1]$, $\mu$ be the Lebesgue measure on $\mathcal B(H)$ (the Borel $\sigma$-algebra) and for any integer $p\ge 1$ let $L^p(H)$ denote the linear
space of all real valued functions $u$ such that
$u^p$ is integrable on $(\mathcal B(H),\mu)$.
To build a contraction map $T$ on $L^2(H)$ we
need $N$-map contractive IFS i.e. a set of maps
$\mathcal W=\{w_1, w_2, \ldots, w_n\}$ and a set of
functions (grey level maps) $\phi=\{\phi_1, \phi_2, \ldots, \phi_n\}$
with $\phi_i:\IR\to\IR$. The operator $T$ corresponding to the
$N$ map IFSM($w$,$\phi$) is
\begin{equation}
(Tu)(x)=\sum_{k=1}^N {'}\phi_k(u(w_k^{-1}(x)))
\label{eq1}
\end{equation}
where the prime means that the sum operates on all those terms
for which $w_k^{-1}$ (the inverse function of $w_k$) is defined.
Let us define the following two sets
$$
{\mathcal Sim}(H)=\{w:H\to H: \exists c\in [0,1), |w(x)-w(y)|=c |x-y|,
\forall x,y\in H\}
$$
$$
{\mathcal Lip}(\IR)=\{\phi:\IR\to\IR:\exists K\in [0,\infty), |\phi(t_1)-\phi(t_2)|\le
K|t_1-t_2|, \forall t_1,t_2\in \IR\}
$$

\begin{theorem}
\cite{Forte}
Let $(w,\phi)$ be an IFSM such that $w_k\in {\mathcal Sim}(H)$ and
$\phi_k\in {\mathcal Lip}(\IR)$ for $1\le K\le N$. Then $T:L^2(H)\to L^2(H)$
and for any $u,v\in L^2(H)$ we have
$$
d_2(Tu,Tv)\le C d_2(u,v)
$$
where
$$
C=\sum_{k=1}^N c_k^{\frac12} K_k.
$$
\end{theorem}

Given $u\in L^2(H)$ the inverse problem consists of finding
the operator $T$ such that
$$
u(x)=(Tu)(x)=
\sum_{k=1}^N {'}\phi_k(u(w_k^{-1}(x))).
$$
In \cite{Forte} it is proved that this problem
 can be reduced to the determination of grey level maps
$\phi_k$ which minimize the collage distance $\Delta^2$
$$
\Delta^2=\|v-Tv\|=
$$
$$
\int_H \sum_{k=1}^{N} {'}\|\phi_k(v(w_k^{-1}(x)))-v(x)\|d\mu.
$$         
In the special case when
\begin{itemize}
\item $\bigcup_{k=1}^N H_k=\bigcup_{k=1}^N w_i(H)=H$ i.e. the sets
$H_k$ ``tile'' $H$
\item $\mu(w_i(H)\cap w_j(H))=0$ for $i\not=j$
\end{itemize}
we say that the maps $w_k$ are nonoverlapping.
Later, in the applications, we will assume that, for $1\le k\le N$,
\begin{itemize}
\item $w_k(x)=s_k x + a_k$
\item $0 < c_k =  |s_k| < 1$ 
\item $\phi_k(t)=\alpha_kt+\beta_k$, $K_k=|\alpha_k|$
\end{itemize}
The collage distance becomes
$$
\Delta^2 = <v-Tv, v-Tv>=
$$
$$
\sum_{k=1}^N \sum_{l=1}^N <\psi_k,\psi_l>\alpha_k \alpha_l +
2 <\psi_k,\psi_l>\alpha_k \beta_l + <\xi_k,\xi_l>\beta_k \beta_l
$$
$$
-2\sum_{k=1}^N <v,\psi_k>\alpha_k+<v,\xi_k>\beta_k+<v,v>
$$
where
$$
\psi_k(x)=v(w_k^{-1}(x)),\qquad \xi_k(x)=I_{w_k(H)}(x)
$$
$\Delta^2$ is a quadratic form in $\alpha_i$ and $\beta_i$, that is
\begin{equation}
\Delta^2=x^T A x + b^T x + c
\label{eq2}
\end{equation}
where $x=(\alpha_1,\ldots \alpha_k, \beta_1, \ldots, \beta_k)$.
The matrix $A$ is symmetric and
$$
a_{i,j}=<\psi_i,\psi_j>, a_{N+i,N+j}=<\xi_i,\xi_j>
$$
$$
a_{i,N+j}=<\psi_i,\xi_j>, b_i=-2<v,\psi_i>, b_{N+i}=-2<v,\xi_i>
$$
and $c=||v||_2^2$.
As in \cite{Forte} we add an additional constraint
in order to guarantee that the minimum of this quadratic
form exists on a compact subset of feasible parameters $\alpha_i$
and $\beta_i$. The additional constraint is
$$
\sum_{k=1}^N c_k (\alpha_k \|v\|_1 + \beta_k)-\|v\|_1\le 0.
$$
The maps $w_k$ are choosen in an infinite set $\mathcal W$ of fixed
affine contraction maps on $H$ which
has the $\mu$-dense and nonoverlapping property (in the
sense of the following definition); When
 $(\alpha_k,\beta_k)=(0,0)$ the corresponding
$w_k$ is superfluous and the $k$-th term can be dropped from \eqref{eq1}.

\begin{definition}
We say that $\mathcal W$ generates a $\mu$-dense and nonoverlapping family
$\A$ of subsets of $H$ if for every $\epsilon>0$ and every $B\subset H$
there exists a finite set of integers $i_k$, $i_k\ge 1$, $1\le k\le N$,
such that

\begin{itemize}

\item $A=\cup_{k=1}^N w_{i_k}(H)\subset B$
\item $\mu(B\backslash A)<\epsilon$
\item $\mu(w_{i_k}(H)\cap w_{i_l}(H))=0$ if $k\not=l$

\end{itemize}
\end{definition}
Let
$$
\mathcal W^N=\{w_1,\ldots w_N\}
$$
be the $N$ truncations of $w$. Let $\Phi^N=\{\phi_1,\ldots,\phi_N\}$
the $N$ vector of affine grey level maps.
Let $x^N$ be the solution of the previous quadratic optimization
problem \eqref{eq2} and $\Delta^2_{N,min}=\Delta^2_N(x^N)$. It can be shown
that $\Delta^2_{N,min}$ may be arbitrarly small when $N\to \infty$ (see \cite{Forte}).

\section{IFSM for stochastic processes on $L^2(H)$}\label{sec:l2}
Let $(\Omega,\A,P)$ be a probability space
and $\{\A_t, t\in H\}$ be a sequence of $\sigma$-algebras such
that $\A_t\subset \A$.
Let $X(\omega,t):\Omega\times H\to \IR$ be a
stochastic process in $L^2(H)$, that is a sequence of
random variables $\A_t$-adapted (that is
each variable $X(\omega,t)$ is $\A_t$-measurable).
Given $\omega\in\Omega$ a trajectory of the
process is the function $X(\omega,t):H\to \IR$ belonging to  $L^2(H)$.
For a given $X(\omega,t)$, the trajectory of the stochastic process, the aim
of the inverse
problem consists in finding the parameters of the IFSM such that
$X(\omega,t)$ is the solution of the equation
$$
X(\omega,t)=T X(\omega,t)\quad\text{for a.a.}\quad \omega\in \Omega
$$
In this case the coefficients
of the matrix $A$ and the vector $b$ of the previous section become
$$
\begin{aligned}
a_{i,j}(\omega)&=\int_H X(\omega,w_i^{-1}(t))X(\omega,w_j^{-1}(t)dt\\
&=
\int_{w_i(H)\cap w_j(H)} X(\omega,w_i^{-1}(t))X(\omega,w_j^{-1}(t)dt
\end{aligned}
$$
and if $i=j$ it becomes
$$
a_{i,i}(\omega)=c_i \int_H X^2(\omega,t)dt
$$
The other elements in the matrix $A$ can be calculated as
$$
a_{N+i,N+j}=<\xi_i,\xi_j>=\int_H I_{w_i(H)}(t)I_{w_j(H)}(t)dt
=\mu(w_i(H)\cap w_j(H))
$$
and
$$
a_{i,N+j}(\omega)=<\psi_i,\xi_j>=\int_{w_i(H)\cap w_j(H)} X(\omega,w_i^{-1}(t))dt
$$
For the vector $b$
$$
b_i(\omega)=-2<X,\psi_i>=\int_H X(\omega,t)X(\omega,w_i^{-1}(t))dt
$$
and
$$
b_{N+i}(\omega)=-2<X,\xi_i>=\int_{w_i(H)} X(\omega,t)dt
$$
In the nonoverlapping case, we have

\begin{itemize}
\item $a_{i,j}=0$, $i\not=j$, and $a_{i,i}(\omega)=c_i \int_H X^2(\omega,t)dt$,
      $1\le i,j\le N$
\item $a_{N+i,N+j}=0$, $1\le i,j\le N$, $i\not=j$ and $a_{N+i,N+i}
=\mu(w_i(H))$
\item $a_{i,N+j}=0$, $1\le i,j\le N$, $i\not=j$ and $a_{i,N+i}(\omega)=c_i
\int_H X(\omega,t)dt$
\end{itemize}
It also holds this self-similarity property.

\begin{theorem}
Let $(\alpha_k,\beta_k)$ be the solution of
the inverse problem with a set of nonverlapping
maps $w_k$ and suppose that $\tilde X(\omega,t)=T \tilde X(\omega,t)$.
Then
$$
\tilde X(\omega,w_i (t+h))-\tilde X(\omega,w_i(t))=
\alpha_i(\tilde X(\omega,t+h)- \tilde X(\omega,t)).
$$
for all $1\le i\le N$.
\end{theorem}

\begin{proof}
In fact we have
$$
\tilde X(\omega,w_i(t+h))-\tilde X(\omega,w_i(t))=
 T \tilde X(\omega,w_i(t+h)) -
T \tilde X(\omega,w_i(t))
$$
$$
\begin{aligned}
&= \sum_{k=1}^{N}
\alpha_k(\tilde X(\omega,w_k^{-1}(w_i(t+h)))+\beta_k - \sum_{k=1}^{N}
\alpha_k(\tilde X(\omega,w_k^{-1}(w_i(t)))+\beta_k \\
&= \alpha_i (\tilde X(\omega,w_i^{-1}(w_i(t+h)))+\beta_i - 
 \alpha_i(\tilde X(\omega,w_i^{-1}(w_i(t)))+\beta_i \\
&=
\alpha_i (\tilde X(\omega,t+h)-\tilde X(\omega,t)).
\end{aligned}
$$
\end{proof}

\subsection{The Kac-Siegert decomposition of $L^2(H)$ stochastic\\ processes}
We suppose that a.e. $X(\omega,t)$ is an element of a
subspace $S$ of $L^2(H)$ and  that $X(\omega,t)$ is a zero-mean process. Let $K$ be the covariance function of this process
that is
$$
K(s,t)=\Cov[X(\omega,s),X(\omega,t)].
$$
and assume
$$
\int_H K(t,t)dt<\infty.
$$
If $\lambda_1\ge \lambda_2\ge ...>0$ comprises
the entire spectrum of eigenvalues of $K$,
where
$$
\int_H f(s)K(s,t)dt=\lambda f(t), \ 0\le t\le 1
$$
and the associated orthonormal eigenfunctions $f_i$ form
a complete set of the subspace $S$ then
the Kac and Siegert decomposition holds:
$$
K(s,t)=\sum_{j=1}^{\infty} \lambda_j f_j(s)f_j(t), \ 0<s,t<1
$$
We also have that
$$
Z_j=\int_H X(t)f_j(t)dt
$$
are uncorrelated random variables with mean $0$ and variance
$\lambda_j$.
The following theorem states some properties of this
decomposition.
 
\begin{theorem}[see Ch.5, \cite{kac}]
Suppose that $X$ and $K$ satisfies the properties above. Then

\begin{itemize}
\item{i)} $\sum_{j=0}^\infty \lambda_j<\infty$
\item{ii)} $\sum_{j=1}^m Z_j f_j\to_{\|\|_2} X$ as $m\to \infty$ a.s. 
\item{iii)} $Z_j=<X,f_j>$ are with mean $0$ and variance $\lambda_j$
\item{iv)} $\int_H X^2(t)dt=\sum_{j=1}^\infty Z_j^2=
\sum_{j=1}^\infty \lambda_j Z_j^{*2}$ where
$Z_j^*=\frac{Z_j}{\sqrt{\lambda_j}}$
\item{v)} $E[X(t)-\sum_{j=1}^m Z_j f_j(t)]^2\to 0$ for each t as $m\to \infty$
\item{vi)} $X=\sum_{j=1}^m \sqrt{\lambda_j}f_j Z_j^*$ with
$Z_j^*=\frac{Z_j}{\sqrt{\lambda_j}}$ uncorrelated with mean $0$ and variance
$\lambda_j$.
\end{itemize}
\end{theorem}

\section{Simulation of Brownian motion via IFSM on $H=[0,1]$}\label{sec:bm}
In the literature there several methods of simulation of the trajectory of the Browian motion, i.e. the stochastic process $\{ B(\omega,t), t\in[0,1]\}$, such that $B(0) = 0$ a.s., $B(t) - B(s)$ is distributed with Gaussian law with zero mean and variance $t-s$, and with independent increments.
\subsubsection*{The Euler method}
In this case, the trajectory is obtained simulating the increments of $B$ in the following way:
$B(0)=0$, $B(t_{i+1}) = B(t_i) + \sqrt{t_{i+1}-t_i} \cdot  Z_i$, where the $Z_i$'s are independent $\mathcal N(0,1)$ random variables. In the other points the trajectory is built by linear
interpolation of these simulated data.

\subsubsection*{The Kac-Siegert method}
Karhunen-Lo\`eve / Kac-Siegert decomposition of $B$ is better for pathwise simulation
$$
B(\omega,t) = \sum_{i=0}^\infty  Z_i \phi_i(t),\quad 0\leq t \leq 1
$$
with
$$
\phi_i(t) = \frac{2\sqrt{2}}{(2i+1)\pi} \sin \left(\frac{(2i+1)\pi t}{2}\right)
$$
$\phi_i$ a basis of orthogonal functions and $Z_i$'s are $\mathcal N(0,1)$
\medskip

\begin{figure}[t]
\begin{center}
\includegraphics{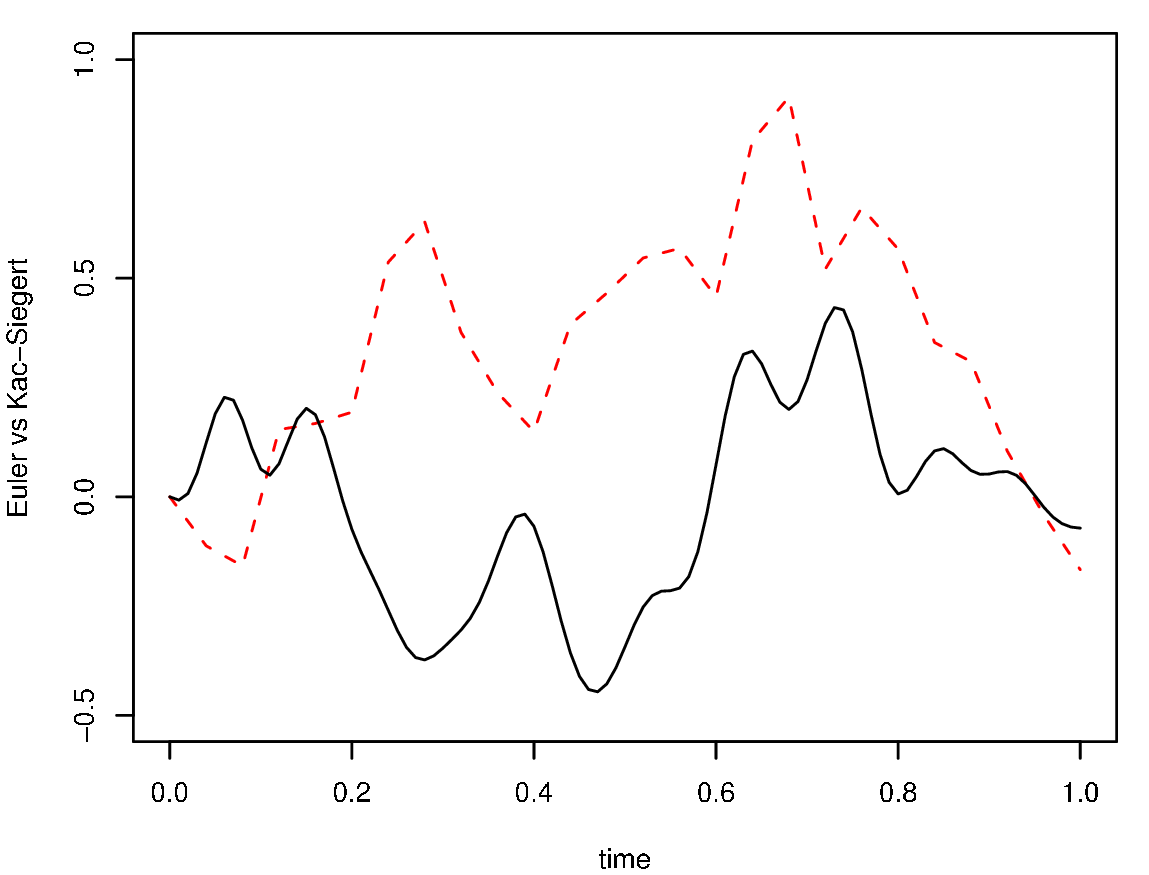}
\end{center}
\caption{Paths of Brownian motion simulated by the Euler scheme (dotted line) and using Kac-Siegert decomposition (continuous line). The  same  ($n=25$) pseudo-random Gaussian  numbers were used.}
\label{fig1}
\end{figure}

The trajectory generated by Euler method is too simple and regular to mimic
the roughness of the BM; moreover the simulated path 
is stochastically equivalent to the true trajectory only on the points of the grid used in the simulation. 
The Kac-Siegert decomposition of the BM is a pathwise approximation which can lead to a too smooth path (see figure \ref{fig1});  our idea is to use IFSM for generating fractal trajectories of the BM.
There are applications in finance (for instance pricing of american options) in which the whole path matters; our IFSM approch produces a global approximation of the trajectory preserving the 
geometric fractal nature of the target. 

This method can be also  used to simulate paths of solutions of stochastic differential equations driven by Brownian motion (e.g. diffusion processes) replacing the linear behaviour  of the Euler trajectory with a fractal object.
For the solution of the inverse problem for the BM, we choose 
the so-called wavelet type maps \cite{Forte}, that is:
$$w_{ij}^*(x)=\frac{x+j-1}{2^i}$$
with $i=1,2,\ldots$ and $j=1\ldots 2^i$
For each fixed $i$, the family $w^*_{ij}$ is a set of nonoverlapping
maps.  For these maps $c_i=2^{-i}<1$. We organize them as follows
$$
w_1 = w^*_{11}\quad w_2 = w^*_{12} \quad w_3 = w^*_{21} \quad w_4 = w^*_{22} \ldots
$$
To simulate a trajectory of $B$ with non overlapping maps we then need to simulate the joint distribution of all this objects
\begin{enumerate}
\item $\displaystyle\int_0^1 B^2(t)dt$
\item $\displaystyle\int_0^1 B(t)dt$
\item $\displaystyle \int_0^1 B(t)B(w_i^{-1}(t))dt = \int_0^1 B(t)B\left(\frac{t-a_i}{s_i}\right)dt  $
\item $\displaystyle\int_{w_i([0,1])} B(t)dt$
\item $\displaystyle\int_0^1 |B(t)| dt$
\end{enumerate}
but it appears to be still a too difficult problem.

\begin{figure}
\begin{center}
\includegraphics[height=.42\textheight,width=\textwidth]{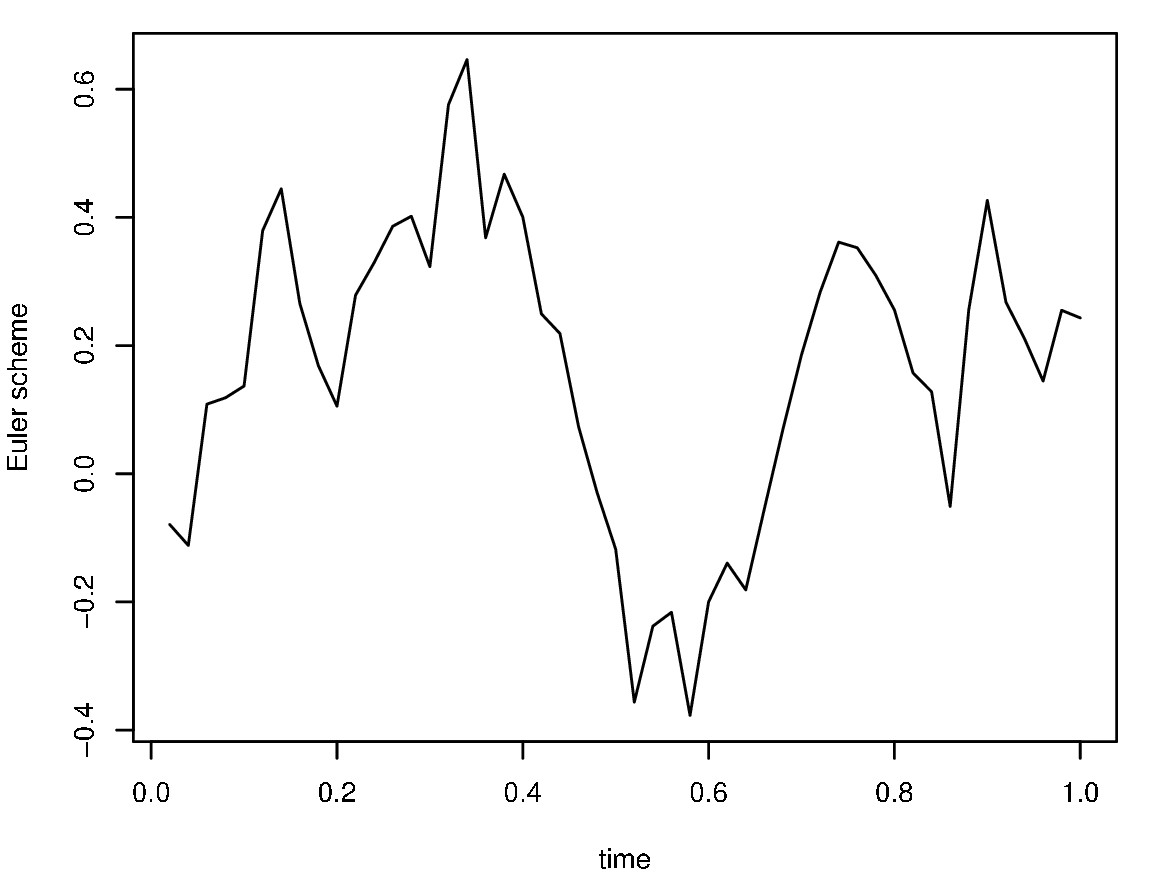}\\
\includegraphics[height=.42\textheight,width=\textwidth]{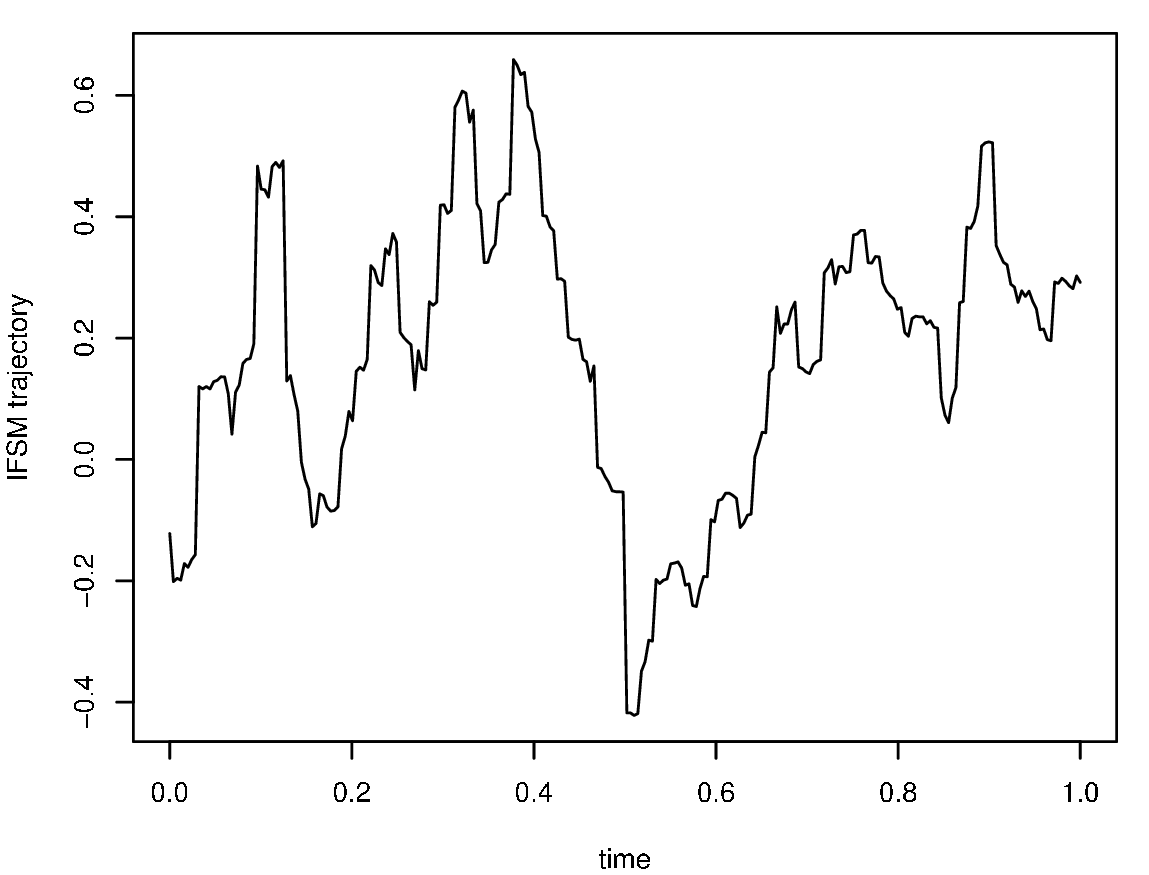}
\end{center}
\caption{Euler scheme versus IFSM trajectory of the Browninan motion. IFSM with wavelet type maps and $M=8$. Both methods used the same 50 Gaussian random terms to generate the trajectory.}
\label{fig2}
\end{figure}

In practice, it is preferable to use all the above maps and not only the subset of non-overlapping maps. In this case, we need simulate the value of the trajectory of the Brownian motion on a fixed grid (using one of the known methods) and we use these points to approximate the integrals in the quadratic form. We then solve the constrained quadratic programming problem using standard algorithms (see e.g. \cite{byrd}).
Figure \ref{fig2} bottom shows an example of trajectory generated using the IFSM approach using wavelet type maps for $i=1, \ldots, M$, $M=8$. Figure \ref{fig2} top represents the Euler trajectory built on 50 Gaussian terms which has been used to build the IFSM. As one can notice the IFSM path shows more ``fractal'' complexity then the corresponding Euler path.

\section{Conclusions}
We have proposed a new method to generate paths of the Brownian motion. These IFSM paths seems to mimic more closely the fractal nature of the trajectory of the Brownian motion then existing schemes.
At current stage we are not able to show formal property of the IFSM path in terms of strong and weak approximation (see \cite{KPS00}). Open source software for generating IFSM trajectories written in C and R language \cite{R} is available via \texttt{ifs} package at {\tt http://CRAN.R-project.org} for free download.

\begin{thebibliography}{99}

\bibitem{bd} Barnsley,
M.F., Demko, S. (1985),  ``Iterated function systems and the
global construction of fractals'', {\sl Proc. Roy. Soc.
London, Ser A}, {\bf 399}, 243-275.

\bibitem{byrd} Byrd, R. H.,
 Lu, P., Nocedal, J. and Zhu, C. (1995),  ``A limited memory
 algorithm for bound constrained optimization'', {\sl SIAM J. Scientific
Computing}, 16, 1190-1208.

\bibitem{Forte} Forte, B.,
 Vrscay, E.R. (1995),  ``Solving the inverse problem for
 function/image approximation using iterated function systems,
 I. Theoretical basis'', {\sl Fractal}, {\bf 2}, 3, 325-334.

\bibitem{fv98} Forte, B., Vrscay,
E.R. (1998),  ``Inverse problem methods for generalized fractal transforms'', in
{\sl Fractal Image Encoding and Analysis}, NATO ASI Series F,
 Vol. 159, ed. Y. Fisher, Springer Verlag, Heidelberg.

\bibitem{hutch} Hutchinson, J. (1981),  ``Fractals
 and self-similarity'', {\sl Indiana Univ.
J. Math.}, {\bf 30}, 5, 713-747.

\bibitem{KPS00} Kloden, P., Platen, E., Shurtz, H. (2000),  Numerical Solution of SDE through computer experiments, Springer, Berlin.

\bibitem{kac} Shorack, G., Wellner, J.A. (1986), Empirical processes with applications to statistics, Wiley, New York.

\bibitem{R}  R Development Core Team (2005), R: A language and environment for
  statistical computing. R Foundation for Statistical Computing,
  Vienna, Austria. ISBN 3-900051-07-0, URL {\tt http://www.R-project.org}

\end{thebibliography}
\end{document}